\documentclass[preprint, 10pt]{elsarticle}
\usepackage{stmaryrd}
\usepackage{amsmath}
\usepackage{tikz,pgfplots}
\usepackage[colorlinks=true,citecolor=blue,linkcolor=blue,urlcolor=blue]{hyperref}
\usepackage[top=1.2in,bottom=1.2in,left=1in, right=1in]{geometry}
\usepackage{amssymb}

\newcommand{\grad}{{\triangledown}}
\newcommand{\ff}{{\mathbf{f}}}
\newcommand{\nn}{{\mathbf{n}}}
\newcommand{\rr}{{\mathbf{r}}}
\newcommand{\uu}{{\mathbf{u}}}
\newcommand{\vv}{{\mathbf{v}}}
\newcommand{\xx}{{\mathbf{x}}}
\newcommand{\yy}{{\mathbf{y}}}
\newcommand{\down}{{\mathrm{down}}}
\newcommand{\exx}{\mathbf{e}_{\xx}}
\newcommand{\esigma}{e_{\sigma}}
\newcommand{\new}{{\mathrm{new}}}
\newcommand{\opt}{{\mathrm{opt}}}
\newcommand{\tsigma}{{\tilde{\sigma}}}
\newcommand{\tvv}{{\tilde{\vv}}}
\newcommand{\txx}{{\tilde{\xx}}}
\newcommand{\up}{{\mathrm{up}}}
\newcommand{\BB}{{\mathcal{B}}}
\newcommand{\DD}{{\mathcal{D}}}
\newcommand{\Div}{{\mathrm{Div}}}
\newcommand{\RR}{{\mathbb{R}}}
\renewcommand{\SS}{{\mathcal{S}}}
\newcommand{\TT}{{\mathcal{T}}}
\newcommand{\scale}{{\mathrm{scale}}}
\newcommand{\sdc}{{\mathrm{sdc}}}

\def\gap{\hspace*{.2in}}

\begin{document}

\begin{frontmatter}

\title{High-order adaptive time stepping for vesicle suspensions with
viscosity contrast}

\author[a]{Bryan Quaife\corref{cor1}} 
\author[a]{George Biros}

\address[a]{Institute for Computational Engineering and Sciences,
University of Texas, Austin, Texas, USA}

\begin{abstract}
We construct a high-order adaptive time stepping scheme for vesicle
suspensions with viscosity contrast.  The high-order accuracy is
achieved using a spectral deferred correction (SDC) method, and
adaptivity is achieved by estimating the local truncation error with
the numerical error of physically constant values.  Numerical examples
demonstrate that our method can handle suspensions with vesicles that
are tumbling, tank-treading, or both.  Moreover, we demonstrate that a
user-prescribed tolerance can be automatically achieved for simulations
with long time horizons.
\end{abstract}

\cortext[cor1]{Corresponding author. Tel.: +1-512-232-3509 ; fax: +1-512-471-3312.}
\end{frontmatter}

\ead{quaife@ices.utexas.edu}

\section{Introduction}
\label{s:Introduction}
Vesicles are deformable capsules filled with an incompressible viscous
fluid. Their mechanical properties are characterized by bending
resistance and tension that enforces local inextensibility of their
membrane.  We are interested in the fluid mechanics of vesicles
suspensions where the bulk fluid is a Newtonian incompressible fluid.
The hydrodynamics of vesicles play an important role in many
biological phenomena~\cite{kra:win:sei:lip1996,sei1997}.  For example,
they are used experimentally to understand properties of
biomembranes~\cite{sac1996} and red blood
cells~\cite{nog:gom2005,poz1990,ghi:rah:bir:mis2011,kao:tah:bib:ezz:ben:bir:mis2011,mis2006}.
Here we extend our work on high-order adaptive time
integrators~\cite{qua:bir2014c} to vesicle suspensions with viscosity
contrast; that is, the viscosity in the interior of each vesicle is
constant, but can differ from the constant viscosity of the exterior
fluid.  This extension is important for applications such as, for
example, simulating microfluidic devices that sort red blood cells of
different viscosity contrasts.  By using an adaptive high-order time
integrator, we can efficiently simulate suspensions to a
user-specified tolerance without requiring a trial and error
procedure.  Moreover, by using a high-order method, fewer time steps
are required when the desired tolerance is small.  This paper is a key
step towards a black-box solver for vesicle suspensions with long time
horizons and that exhibit dynamics of varying complexity.

The most significant limitation of this work is that it is implemented
only in two dimensions. 
Fortunately, the time integrator and adaptivity strategy we introduce is
independent of the dimension and we plan to extend this work to
three dimensions in the future.  While we now allow for variable time
step sizes, we do not address multirate time integrators.  That is, at
each time step, the same time step size is applied to all the vesicles.
In addition to time adaptivity, spatial adaptivity also needs to be
addressed in the future.  Since vesicles may develop and smooth regions
of high curvature, allowing for only one spatial grid is inefficient.
In this paper, the spatial domain is fixed for the entire simulation
but, in the numerical examples, when performing a convergence study in
time, we choose a spatial discretization that is sufficiently large that
the temporal error dominates.

Vesicle suspensions, in two and three dimensions, have been well-studied
and we refrain from doing a thorough literature review of their physics.
Instead, we focus on work related to our time integrator.  The spectral
deferred correction (SDC) method was first introduced but Dutt,
Greengard, and Rokhlin~\cite{dut:gre:rok2000}.  In its original form, it
was a stable method to construct high-order solutions of initial value
problems.  It was extended to implicit-explicit (IMEX) methods by
Minion~\cite{min2003}, where the splitting between stiff and non-stiff
terms is additive.  SDC has also been used as a parallel-in-time time
stepper~\cite{chr:mac:ong:spi2014}, but this work only considers initial
value problems which are much cheaper to solve than the
integro-differential equations we solve.  SDC has been used for an
integro-differential equation~\cite{hua:lai:xia2006} that simulates a
diffusion process with moving interfaces, but their governing equations
are less stiff and the interface equations are simpler than ours.  To
the best of our knowledge, our work~\cite{qua:bir2014c} is the first to
apply SDC to an operator-based IMEX splitting of the equations governing
vesicle suspensions.  It is this work that we now extend to vesicles
with viscosity contrast.

The remainder of the paper is organized as follows.  In
Section~\ref{s:Formulation} we present equations that govern vesicle
dynamics, reformulate them in terms of integral equations, and introduce
SDC.  In Section~\ref{s:NumericalScheme}, the numerical schemes for the
governing equations are presented and we discuss the adaptive time
stepping strategy.  In Section~\ref{s:NumericalExamples}, we present two
numerical results, and we make concluding remarks in
Section~\ref{s:Conclusion}.

\section {Formulation}
\label{s:Formulation}
Neglecting inertial forces, the dynamics of a vesicle is fully
characterized by the position of the interface $\xx(s,t) \in \gamma$,
where $s$ is the arclength, $t$ is time, and $\gamma$ is the membrane
of the vesicle.  The position is determined by solving a moving
interface problem that models the mechanical interactions between the
viscous incompressible fluid in the exterior and interior of the
vesicle and the vesicle membrane, all the while, requiring that the
membrane maintains its length ({\em inextensibility condition}).  In
addition to the position $\xx(s,t)$, the other main variables are the
fluid velocity $\uu$, the fluid stress $T$, the pressure $p$, the
membrane tension $\sigma$, and the stress jump across $\gamma$,
$\ff(\xx) = \llbracket T \rrbracket \mathbf{n}$, where $\mathbf{n}$ is
the unit outward normal.  The stress jump is equal to the sum of a
force due to the vesicle membrane bending modulus $\kappa_{b}$ and a
force due to the tension $\sigma$.

Given these definitions, the equations for an unbounded suspension of
vesicles are
\begin{align}
\label{e:vesicles-pde}
  \begin{alignedat}{2}
    \mu \grad \cdot (\grad \uu + \grad \uu^{T}) = \grad p(\xx), 
    &\hspace{20pt} \xx \in \RR^{2}\backslash\gamma, 
      \gap &&  \mbox{conservation of momentum}, \\
    \grad \cdot \uu(\xx) = 0, &\hspace{20pt} 
    \xx \in \RR^{2} \backslash \gamma, 
      \gap &&\mbox{conservation of mass}, \\
    \uu(\xx,t) = \vv_{\infty}(\xx), &\hspace{20pt} |\xx| 
      \rightarrow \infty, 
       \gap  &&\mbox{far-field velocity}, \\
    \dot{\xx}(t) = \uu(\xx,t), &\hspace{20pt} \xx \in \gamma,
      \gap   &&\mbox{velocity continuity},\\
    \llbracket T \rrbracket \mathbf{n} = -\kappa_{b}\xx_{ssss} 
    + (\sigma(\xx) \xx_{s})_{s}, &\hspace{20pt} \xx \in \gamma,
      \gap  &&\mbox{nonzero stress jump}, \\
    \xx_{s} \cdot \uu_{s} = 0, &\hspace{20pt} \xx \in \gamma,
      \gap &&\mbox{membrane inextensibility}. \\
  \end{alignedat}
\end{align}
While the viscosity is taken to be constant inside each vesicle, these
values can differ from the viscosity of the exterior fluid.  In
particular, we define the viscosity contrast $\nu_{p} =
\mu_{p}/\mu_{0}$.

Since the viscosity is piecewise constant with a discontinuity along the
interface, we use an integral equation formulation using the Stokes
free-space Green's function~\cite{rah:vee:bir2010}.  We introduce the
following integral and differential operators that we use to
formulate~\eqref{e:vesicles-pde} in terms of unknowns that are defined
only on the vesicle boundary $\gamma$:
\begin{itemize}[]
  \item{\em Single- and double-layer potentials:} 
  \begin{align*} 
    \SS(\gamma_{j},\gamma_{k}) \ff &:= \frac{1}{4\pi\mu_{0}}\int_{\gamma_{k}}\left(
      -\boldsymbol{I} \log \rho
    + \frac{\rr \otimes \rr}{\rho^{2}} \right)\ff\, ds_{\yy},
    && \xx \in \gamma_{j}, \\
    \DD(\gamma_{j},\gamma_{k}) \uu &:= \frac{1-\nu_{k}}{\pi}\int_{\gamma_{k}}
      \frac{\rr \cdot \nn}{\rho^{2}}\frac{\rr \otimes \rr}{\rho^{2}}\uu\,
    ds_{\yy}, && \xx \in \gamma_{j},
  \end{align*}
  where $\rr = \xx - \yy$, $\rho = \|\rr\|$, and $\nn$ is the outward
  unit normal of $\gamma_{k}$.
  \item{\em Bending, tension, and surface divergence:} 
  \begin{align*}
    \BB(\gamma_{k}) \ff := \kappa_{b}\frac{d^{4}}{ds^4} \ff, \quad
    \TT(\gamma_{k}) \sigma := (\sigma \xx_{s})_{s}, \quad
    \Div(\gamma_{k}) \ff := \frac{d\xx}{ds} \cdot \frac{d\ff}{ds},
  \end{align*}
  where each arclength derivative is taken with respect to $\xx_{k}$,
  and $\kappa_{b}$ is the bending modulus.
\end{itemize}
Note that all these operators are linear once $\gamma_{j}$ and
$\gamma_{k}$ are fixed.

\subsection{Integral equation formulation}
Let $\omega_{j}$, $j=1,\dots,M$, be the interior of vesicle $j$.  An
integral equation representation of the fluid velocity is~\cite{poz1992}
\begin{align*}
  \alpha \uu(\xx) = \vv_{\infty}(\xx) + \sum_{k=1}^{M}\SS(\xx,\xx_{k})
  \left(-\BB(\xx_{k})\xx_{k} + \TT(\xx_{k})\sigma_{k}\right) + 
  \sum_{k=1}^{M}\DD(\xx,\xx_{k})\uu(\xx_{k}), \quad \xx \in \RR^{2},
\end{align*}
where
\begin{align*}
  \alpha = \left\{
  \begin{array}{ll}
    \nu_{j}             &\gap \xx \in \omega_{j}, \\
    \frac{1+\nu_{j}}{2} &\gap \xx \in \gamma_{j}, \\
    1                   &\gap \xx \in 
          \RR^{2} \backslash \overline{\omega_{j}}.
  \end{array}
  \right.
\end{align*}
Applying the no-slip boundary condition on the boundary of the vesicle,
we have
\begin{align*}
  \alpha_{j}\frac{d\xx_{j}}{dt} = \vv_{\infty}(\xx_{j}) + 
  \sum_{k=1}^{M}\SS(\xx_{j},\xx_{k})
  \left(-\BB(\xx_{k})\xx_{k} + \TT(\xx_{k})\sigma_{k}\right) + 
  \sum_{k=1}^{M}\DD(\xx_{j},\xx_{k})\frac{d\xx_{k}}{dt}, \quad 
  \alpha_{j} = \frac{1+\nu_{j}}{2}.
\end{align*}
The tension $\sigma$ acts as a Lagrange multiplier to impose the
inextensibility condition $\Div(\xx_{j}) \frac{d\xx_{j}}{dt} = 0$.  In
order to ease the presentation, we abuse notation by dropping the
subscripts and summations, and we write
\begin{align}
  \alpha \frac{d\xx}{dt} = \vv_{\infty}(\xx) + 
    \SS(\xx)\left(-\BB(\xx)\xx + \TT(\xx)\sigma\right) + \DD(\xx)\frac{d\xx}{dt},
  \quad \alpha = \frac{1+\nu}{2},
  \label{e:vesicleDynamics}
\end{align}
with the inextensibility condition $\Div(\xx)\frac{d\xx}{dt} = 0$.

In the SDC framework, it is convenient to
reformulate~\eqref{e:vesicleDynamics} as the Picard integral 
\begin{align}
  \xx(t) = \xx_{0} + \int_{0}^{t} (\alpha I - \DD(\xx))^{-1}
    (\vv_{\infty}(\xx)-\SS(\xx)(\BB(\xx)\xx + \TT(\xx)\sigma)) d\tau.
  \label{e:picard}
\end{align}
A similar Picard integral can be derived for confined flows by replacing
$\vv_{\infty}$ with a double-layer potential with an unknown density
function defined on the solid walls~\cite{qua:bir2014c}.  Then, a
no-slip boundary condition on the solid walls results in a second-kind
integral equation that must also be satisfied.  In this work, for
simplicity, we only present results for unbounded flows.

\subsection{Spectral Deferred Correction}
Spectral deferred correction (SDC) is an iterative method for solving
Picard integral equations such as~\eqref{e:picard}.  It was first
introduced as a stable deferred correction method for solving initial
value problems~\cite{dut:gre:rok2000}.  Our SDC formulation is most
closely related to the work of Minion~\cite{min2003} who was the first
to investigate the coupling of SDC with IMEX time integrators.  We first
compute a provisional solution $\txx$ and $\tsigma$ using some time
integrator, and then form the residual of~\eqref{e:picard}
\begin{align}
  \rr(t) := \xx_{0} - \txx(t) + \int_{0}^{t} 
      \tvv(\tau) d\tau,
  \label{e:residual}
\end{align}
where the provisional velocity is $\tvv = (\alpha I -
\DD(\txx))^{-1}(\vv_{\infty}(\txx) - \SS(\txx)\BB(\txx)\txx +
\SS(\txx)\TT(\txx)\tsigma)$.  We now define the errors in position,
$\exx := \xx - \txx$, and tension, $\esigma := \sigma - \tsigma$.  By
substituting $\exx$ and $\esigma$ into~\eqref{e:picard}, we have
\begin{align*}
  \exx(t) = \xx_{0} - \txx(t) + \int_{0}^{t} 
    (\alpha I - \DD(\txx+\exx))^{-1}(\vv_{\infty}(\txx+\exx)
    -\SS(\txx+\exx)\BB(\txx+\exx)(\txx+\exx) 
    +\SS(\txx+\exx)\TT(\txx+\exx)(\tsigma+\esigma)) d\tau.
\end{align*}
Following the usual SDC framework, we introduce the residual into the
error equation
\begin{align}
\begin{split}
  \exx(t) &= \rr(t) + \int_{0}^{t}\left\{
    (\alpha I-\DD(\txx+\exx))^{-1}\vv_{\infty}(\txx+\exx) - 
    (\alpha I-\DD(\txx))^{-1}\vv_{\infty}(\txx)\right\}d\tau \\
  &+\int_{0}^{t}\left\{
    -(\alpha I-\DD(\txx+\exx))^{-1}
    \SS(\txx+\exx)\BB(\txx+\exx) +
    (\alpha I-\DD(\txx))^{-1}\SS(\txx)\BB(\txx)\right\}\txx d\tau \\
  &+\int_{0}^{t}\left\{
    (\alpha I-\DD(\txx+\exx))^{-1}
    \SS(\txx+\exx)\TT(\txx+\exx) -
    (\alpha I-\DD(\txx))^{-1}\SS(\txx)\TT(\txx)\right\}\tsigma d\tau \\
  &+\int_{0}^{t} (\alpha I-\DD(\txx+\exx))^{-1}
      \left\{-\SS(\txx+\exx)\BB(\txx+\exx)\exx + 
      \SS(\txx+\exx)\TT(\txx+\exx)\esigma\right\}d\tau.
\end{split}
\label{e:errorPicard}
\end{align}
Finally, the inextensibility constraint is
$\Div(\txx+\exx)\frac{d(\txx+\exx)}{dt} = 0$.

\section{Numerical Scheme}
\label{s:NumericalScheme}
Several challenges arise when forming numerical solutions
of~\eqref{e:picard} (and~\eqref{e:errorPicard}).  For instance, the
bending operator $\BB$ introduces stiffness which is resolved by using a
semi-implicit discretization of~\eqref{e:picard}.  Another challenge is
to accurately approximate the layer potentials $\SS$ and $\DD$ for
vesicles that are arbitrarily close to one another which is resolved
with specialized quadrature that efficiently handles near-singular
integrals.  Additional complications and more details are discussed in
previous work~\cite{vee:gue:zor:bir2009,qua:bir2014b}.  Here we focus on
the numerical aspects associated with
solving~\eqref{e:picard},~\eqref{e:residual},~\eqref{e:errorPicard}, and
with computing optimal time step sizes.

SDC is an iterative method that computes successively more accurate
solutions of~\eqref{e:picard}.  The SDC iteration is
\begin{enumerate}
  \item Form a provisional solution $\txx$ and $\tsigma$
  of~\eqref{e:picard}
  \item Form the residual~\eqref{e:residual} of the provisional
  solution.
  \label{i:step2}
  \item Compute an approximate solution of the error
  equation~\eqref{e:errorPicard}.
  \item Define the provisional solution as $\txx + \exx$ and
  $\tsigma{} + \esigma$.
  \item If desired, return to step~\ref{i:step2}.
\end{enumerate}
We have numerically shown~\cite{qua:bir2014c} that if~\eqref{e:picard}
and~\eqref{e:errorPicard} are solved with first-order methods, then
after $n_{\sdc}$ SDC iterations, the order of accuracy of the solution
is $n_{\sdc} + 1$, as long as this order does not exceed the order of
accuracy of the quadrature formula for $\rr$.

\subsection{Discretization of the Picard integral}
A first-order approximation of~\eqref{e:picard} is formed using the
semi-implicit time integrator
\begin{align}
  \alpha\frac{\xx^{N+1} - \xx^{N}}{\Delta t} = \vv_{\infty}^{N} +
  \SS^{N}\left(-\BB^{N}\xx^{N+1} + \TT^{N}\sigma^{N+1}\right) +
  \DD^{N}\left(\frac{\xx^{N+1} - \xx^{N}}{\Delta t}\right).
  \label{e:firstOrderDis}
\end{align}
The time integrator~\eqref{e:firstOrderDis} treats both intra- and
inter-vesicle interactions semi-implicitly, with the latter resolving
stiffness introduced by vesicles that are close to one another.  We
solve~\eqref{e:firstOrderDis} with the generalized minimal residual
method (GMRES) with an exact block-diagonal preconditioner.  This method
has been thoroughly tested for stability, first-order convergence, and
GMRES iteration mesh-independence on a variety of vesicle
suspensions~\cite{qua:bir2014b}.

To compute the residual~\eqref{e:residual}, we
approximate~\eqref{e:picard} at $p$ substeps of $[0,\Delta t]$, and
then use an appropriate quadrature formula to compute $\rr$.  To
achieve high-order accuracy, avoid extrapolation, and minimize the
Runge phenomenon, we use Gauss-Lobatto points, whose quadrature error
is $\mathcal{O}(\Delta t^{2p-3})$.  Once a provisional first-order
solution and the residual are computed, the error
equation~\eqref{e:errorPicard} is numerically solved for $\exx$ and
$\esigma$.  We discretize~\eqref{e:errorPicard} by dropping all the
terms except the first and the last
\begin{align}
  \alpha\frac{\exx^{N+1} - \exx^{N}}{\Delta t} = 
  \left(\alpha I - \DD^{N+1}\right)
  \left(\frac{\rr^{N+1} - \rr^{N}}{\Delta t}\right) +
  \SS^{N+1}\left(-\BB^{N+1}\exx^{N+1} + \TT^{N+1}\esigma^{N+1}\right) +
  \DD^{N+1}\left(\frac{\exx^{N+1} - \exx^{N}}{\Delta t}\right),
  \label{e:firstOrderCor}
\end{align}
and we have observed numerically~\cite{qua:bir2014c} that each solution
of~\eqref{e:firstOrderCor} results in an additional order of
convergence.  In~\eqref{e:firstOrderCor}, the operators are discretized
at $\txx^{N+1}$ which is computed when the provisional solution is
formed.  Now, equations~\eqref{e:firstOrderDis}
and~\eqref{e:firstOrderCor} require inverting the same linear system;
the only differences are the right hand side and the state at which the
operators are discretized.  We also note that if this iteration
converges, that is, $\exx \rightarrow 0$ and $\esigma \rightarrow 0$,
then $\rr^{N+1} = \rr^{N}$, implying that the Picard
integral~\eqref{e:picard} has been solved up to quadrature error.

\subsection{Discretization of the inextensibility condition}
Based on our first-order discretization~\eqref{e:firstOrderDis}, the
inextensibility condition $\Div(\xx)(d\xx/dt) = 0$ for the provisional
solution is equivalent to $\Div(\xx^{N}) \xx^{N+1} = 1$.  Once the
provisional solution $\txx$ and $\tsigma$ have been computed, the
inextensibility condition is $\Div(\txx+\exx)\frac{d(\txx+\exx)}{dt} =
0$, which, using the first-order discretization~\eqref{e:firstOrderCor}
is equivalent to
\begin{align*}
  \Div^{N+1} (\exx^{N+1}) = \Div^{N+1} (\exx^{N}) + 
  \Div^{N+1} (\rr^{N+1} - \rr^{N} + \Delta t_{N} \tvv^{N+1}),
\end{align*}
where $\Delta t_{N}$ is the size of the substep.  Unfortunately, we
have experimentally found that this discretization requires a very
small time step before SDC converges.  As an alternative, we formulate
the inextensibility constraint as
\begin{align*}
  (\txx^{N+1} + \exx^{N+1})_{s_{0}} \cdot 
  (\txx^{N+1} + \exx^{N+1})_{s_{0}}  = 1,
\end{align*}
where $s_{0}$ is the arclength component of the vesicle at the initial
Gauss-Lobatto point.  Assuming that the error $\exx$ is much smaller
than $\txx$, we can neglect the quadratic term and enforce
inextensibility by requiring
\begin{align}
  \txx^{N+1}_{s_{0}} \cdot (\exx^{N+1})_{s_{0}} = \frac{1}{2}
    (1 - \txx_{s_{0}}^{N+1} \cdot \txx_{s_{0}}^{N+1}).
  \label{e:errorDivergence}
\end{align}
Because the arclength term is discretized at $s_{0}$, the
inextensibility condition for the error~\eqref{e:errorDivergence}
differs slightly from the inextensibility condition for the provisional
solution.  However, the resulting system converges for much larger
values of $\Delta t$ and, if $\exx \rightarrow 0$, then the
inextensibility condition $\txx^{N+1}_{s_{0}} \cdot \txx^{N+1}_{s_{0}} =
1$ is exactly satisfied.

\subsection{Time step size selection}
During the course of a vesicle simulation, it is advantageous to adjust
the time step size based on the complexity of the dynamics.  For
instance, when a vesicle has regions of high curvature or two vesicles
are close to one another, a smaller time step should be taken.  The main
advantages of an adaptive time step size are: the simulation is sped up
since the largest allowable time step size is taken; and, no trial and
error procedure is necessary to find a time step size that achieves a
desired tolerance.

We adopt the common strategy of controlling an estimate of the local
truncation error~\cite{hai:nor:wan1993}.  Because of the
inextensibility and incompressibility conditions, vesicle suspensions
have a natural estimate of the local truncation error: the errors in
area and length.  The advantage of this estimate is that multiple
numerical solutions do not need to be formed to estimate the local
truncation error.  For a single vesicle (multiple vesicles are handled
by taking the maximum over all the errors in area and length), let the
area and length at time $t$ be $A(t)$ and $L(t)$, and suppose that the
desired tolerance for the global error is $\epsilon$ at the normalized
time horizon $T=1$.  We compute the solution at time $t+\Delta t$ using
a time stepping scheme and then compute the new area $A(t + \Delta t)$
and length $L(t + \Delta t)$.  The solution is accepted only if
\begin{align}
  |A(t+\Delta t) - A(t)| \leq \frac{A(t) \Delta t}{1-t} 
  \left(\epsilon - \frac{|A(t) - A(0)|}{A(t)}\right),
  \label{e:errorBounds}
\end{align}
and a similar condition for the length.  Instead of the usual condition
$|A(t+\Delta t) - A(t)| \leq A(t) \Delta t \epsilon$, we are adjusting
the amount of local truncation error that can be committed in
$[t,t+\Delta t]$ based on the error committed in $[0,t]$.  Using this
strategy, our scheme can commit larger local truncation errors than
$\epsilon \Delta t$ if the error committed in $[0,t]$ is less than
$\epsilon t$.

Regardless of the acceptance or rejection of the solution at time $t +
\Delta t$, a new time step size must be chosen.  Based on an asymptotic
argument~\cite{qua:bir2014c}, to make the bound
in~\eqref{e:errorBounds} tight, the optimal time step size is
\begin{align*}
  \Delta t_{\opt} = \left\{
    \frac{A(t)}{|A(t+\Delta t) - A(t)|} \frac{\Delta t}{1-t} \left(
    \epsilon - \frac{|A(t)-A(0)|}{A(t)} \right)
    \right\}^{1/k} \Delta t,
\end{align*}
where $k$ is the order of the time integrator.  A similar optimal time
step size is computed based on the length and the smaller of these
values is used to compute the next time step size.  Since the choice of
$\Delta t_{\opt}$ assumes that the error is in the asymptotic regime,
we impose several safety factors on the new time step size $\Delta
t_{\new}$.  First, we require $\Delta t_{\new} / \Delta t \in
[\beta_{\down},\beta_{\up}]$, where $\beta_{\down} = 0.6$ and
$\beta_{\up} = 1.5$.  That is, the time step size is not allowed to
change too quickly.  Next, we multiply the new time step size by a
safety factor $\beta_{\scale} = \sqrt{0.9}$ to increase the likelihood
of the next time step size being accepted.  These three parameters are
chosen from experience and can be changed if required by a certain
application.  As a final restriction, we never increase the time step
size if the previous time step size is
rejected~\cite{hai:nor:wan1993}.  In summary, if the previous time step
is accepted, the new time step size is
\begin{align*}
  \Delta t_{\new} = \beta_{\scale}^{1/k}\min(\beta_{\up}\Delta t,
    \max(\Delta t_{\opt},\beta_{\down}\Delta t)),
\end{align*}
and if the previous time step is rejected, the new time step size is
\begin{align*}
  \Delta t_{\new} = \beta_{\scale}^{1/k}\min(\Delta t,
    \max(\Delta t_{\opt},\beta_{\down}\Delta t)).
\end{align*}

\section{Numerical Examples}
\label{s:NumericalExamples}
We consider vesicles with reduced area 0.66, initially parameterized
as $\xx(\theta,0) = (\cos \theta,3\sin \theta)$, in the shear flow
$\vv_{\infty} = (y,0)$.  We report the number of uniform time steps,
$m$, the error in area, $e_{A}$, and the error in length, $e_{L}$.  To
evaluate the cost of the algorithm, we report the number of
matrix-vector-multiplications, (``{\em matvecs}'' for short),
required to iteratively solve~\eqref{e:firstOrderDis}
and~\eqref{e:firstOrderCor} and to compute the provisional velocity
$\tvv$, and the CPU time required by Matlab on a six-core 2.67GHz
Intel Xeon processor with 24GB of memory.

We first verify that our time integrator with one SDC correction
provides second-order convergence.  Then, we use this time integrator
with an adaptive time step size for two different problems.  We use
$p=5$ Gauss-Lobatto substeps which achieves $7^{th}$-order accuracy
when computing the residual $\rr(t)$.  We set the block-diagonal
preconditioned GMRES tolerance to 1E-10 which, based on previous
work~\cite{qua:bir2014b,qua:bir2014c}, is the smallest allowable
tolerance due to the ill-conditioning of the governing equations.

\subsection{Tank-treading and tumbling}
We consider a single vesicle with two different viscosity contrasts,
and discretized with $N=96$ points, which is sufficiently refined to
guarantee that the temporal error dominates.  We use a fixed time step
size to verify that one SDC correction results in second-order
convergence.  Table~\ref{t:shear_nu4_convStudy} corresponds to a
tank-treading vesicle ($\nu=4$) while
Table~\ref{t:shear_nu15_convStudy} corresponds to a tumbling vesicle
($\nu=15$).  With the prescribed time horizon, the tumbling vesicle
makes approximately 2.5 rotations.  In both regimes, we achieve in
excess of second-order convergence but, in general, we only expect
second-order convergence.

\begin{table}[h]
\caption{\label{t:shear_nu4_convStudy} The number of uniform time steps,
errors in area and length, the number of {\em matvecs}, and the CPU time
for a single vesicle in the tank-treading regime ($\nu=4$).  We see that
second-order (in fact, it is closer to third-order) convergence is
achieved.}
\begin{tabular*}{\hsize}{@{\extracolsep{\fill}}lllll@{}}
\hline
$m$ & $e_{A}$ & $e_{L}$ & {\em matvecs} & CPU \\
\hline
75   & 1.39E-2 & 1.43E-2 & 9.45E3 & 8.16E1 \\
150  & 1.42E-3 & 1.68E-3 & 1.60E4 & 1.53E2 \\
300  & 1.14E-4 & 2.03E-4 & 2.81E4 & 2.88E2 \\
600  & 8.33E-7 & 2.46E-5 & 5.16E4 & 5.68E2 \\
1200 & 3.80E-6 & 2.99E-6 & 9.56E4 & 1.08E3 \\
\hline
\end{tabular*}
\end{table}

\begin{table}[h]
\caption{\label{t:shear_nu15_convStudy} The number of uniform time
steps, errors in area and length, the number of {\em matvecs}, and the
CPU time for a single vesicle in the tumbling regime ($\nu=15$).
Again, second-order (in fact, it is closer to 2.5-order) convergence is
achieved.  Also, we see that GMRES requires fewer iterations to solve
suspensions with larger viscosity contrasts, but larger global errors
are incurred.}
\begin{tabular*}{\hsize}{@{\extracolsep{\fill}}lllll@{}}
\hline
$m$ & $e_{A}$ & $e_{L}$ & {\em matvecs} & CPU \\
\hline
75   & 1.25E-1 & 3.80E-2 & 8.40E3 & 7.94E1 \\
150  & 2.07E-2 & 5.22E-3 & 1.45E4 & 1.55E2 \\
300  & 3.27E-3 & 6.28E-4 & 2.57E4 & 2.83E2 \\
600  & 5.51E-4 & 7.88E-5 & 4.68E4 & 5.38E2 \\
1200 & 1.04E-4 & 9.89E-6 & 8.73E4 & 1.09E3 \\
\hline
\end{tabular*}
\end{table}

We now test our adaptive time stepper with viscosity contrast $\nu=4$
(Table~\ref{t:shear_nu4_sdc1}), $\nu=10$
(Table~\ref{t:shear_nu10_sdc1}), and $\nu=15$
(Table~\ref{t:shear_nu15_sdc1}).  In all of the examples, without any a
priori knowledge of the time step size, the desired tolerance is
achieved.  Unsurprisingly, larger time steps can be taken when the
vesicle is tank-treading ($\nu=4$) as opposed to tumbling ($\nu=10$ and
$\nu=15$).  In all three simulations, as the tolerance decreases, the
number of required time steps grows unacceptably large.  This can be
resolved by taking two SDC corrections (bottom rows of
Tables~\ref{t:shear_nu4_sdc1}--\ref{t:shear_nu10_sdc1}) at the expense
of solving~\eqref{e:errorPicard} twice instead of once which, in theory,
results in a third-order time integrator.  This additional solve is
justified since the reduction in CPU time is 41\% for $\nu=4$, 60\% for
$\nu=10$, and 63\% for $\nu=15$.  Therefore, we conclude that as more
tolerance is requested, the additional cost of using additional SDC
iterations is justified.  In order to continue our development of a
black-box solver for vesicle suspensions, a systematic method for
choosing an appropriate number of SDC iterations based on the requested
tolerance is under investigation.

\begin{table}[h]
\caption{\label{t:shear_nu4_sdc1}Adaptive time stepping with $\nu = 4$.
The adaptive time step size does a good job of achieving the desired
tolerance with only a few rejected time steps.  The final column used
two, rather than one, SDC corrections.}
\begin{tabular*}{\hsize}{@{\extracolsep{\fill}}lllllll@{}}
\hline
Tolerance & $e_{A}$ & $e_{L}$ & Accepts & Rejects & {\em matvecs} & CPU \\
\hline
1E-1  & 9.93E-2 & 4.27E-2 & 54   & 2  & 7.92E3 & 6.75E1 \\
1E-2  & 9.94E-3 & 2.67E-3 & 131  & 4  & 1.45E4 & 1.40E2 \\
1E-3  & 9.98E-4 & 1.10E-4 & 381  & 11 & 3.51E4 & 3.77E2 \\
1E-4  & 1.00E-4 & 4.04E-6 & 1171 & 17 & 9.28E4 & 1.15E3 \\
1E-4* & 1.00E-4 & 6.65E-7 & 316  & 5  & 4.59E4 & 6.75E2 \\
\hline
\end{tabular*}
\end{table}

\begin{table}[h]
\caption{\label{t:shear_nu10_sdc1}Adaptive time stepping with $\nu =
10$.  While the tolerance is always achieved, the final error is not
tight to the desired error.  This is a result of the error decreasing
near the time horizon (see Fig.~\ref{f:shear1VesErrors}).  The final
column used two, rather than one, SDC corrections.}
\begin{tabular*}{\hsize}{@{\extracolsep{\fill}}lllllll@{}}
\hline
Tolerance & $e_{A}$ & $e_{L}$ & Accepts & Rejects & {\em matvecs} & CPU \\
\hline
1E-1  & 4.27E-2 & 1.01E-2 & 86   & 36  & 1.39E4 & 1.31E2 \\
1E-2  & 3.62E-3 & 1.10E-3 & 230  & 92  & 2.93E4 & 3.19E2 \\
1E-3  & 2.67E-4 & 1.29E-4 & 678  & 169 & 6.64E4 & 7.76E2 \\
1E-4  & 3.61E-5 & 2.31E-5 & 2079 & 100 & 1.51E5 & 2.16E3 \\
1E-4* & 5.87E-5 & 1.77E-6 & 310  & 120 & 5.77E4 & 8.73E2 \\
\hline
\end{tabular*}
\end{table}

\begin{table}[h]
\caption{\label{t:shear_nu15_sdc1}Adaptive time stepping with $\nu =
15$.  The adaptive time step size does a good job of achieving the
desired tolerance, but because of the more complicated dynamics of
tumbling, there are more rejected time steps when compared to the
tank-treading vesicle.  The final column used two, rather than one, SDC
corrections.}
\begin{tabular*}{\hsize}{@{\extracolsep{\fill}}lllllll@{}}
\hline
Tolerance & $e_{A}$ & $e_{L}$ & Accepts & Rejects & {\em matvecs} & CPU \\
\hline
1E-1  & 9.71E-2 & 6.22E-3 & 91   & 43  & 1.49E4 & 1.41E2 \\
1E-2  & 9.65E-3 & 7.22E-4 & 245  & 95  & 3.02E4 & 3.28E2 \\
1E-3  & 9.77E-4 & 5.50E-5 & 725  & 186 & 7.01E4 & 8.47E2 \\
1E-4  & 9.82E-5 & 3.45E-6 & 2227 & 142 & 1.63E5 & 2.35E3 \\
1E-4* & 6.13E-5 & 5.11E-7 & 333  & 123 & 5.99E4 & 9.37E2 \\
\hline
\end{tabular*}
\end{table}

\begin{figure}[htps]
\centering
\begin{minipage}{.48\textwidth}
\includegraphics{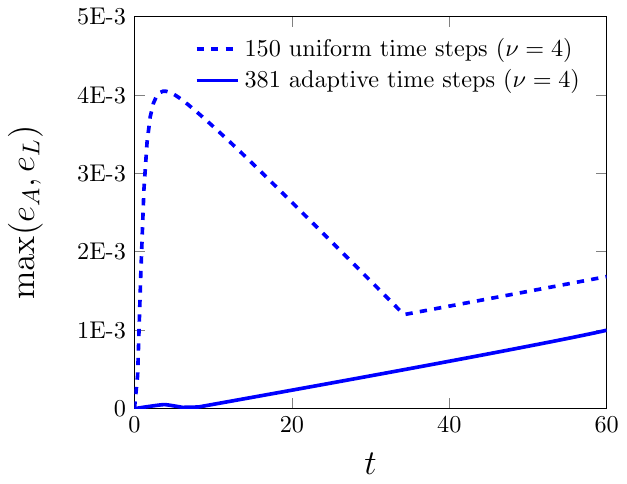}
\end{minipage}
\begin{minipage}{.48\textwidth}
\includegraphics{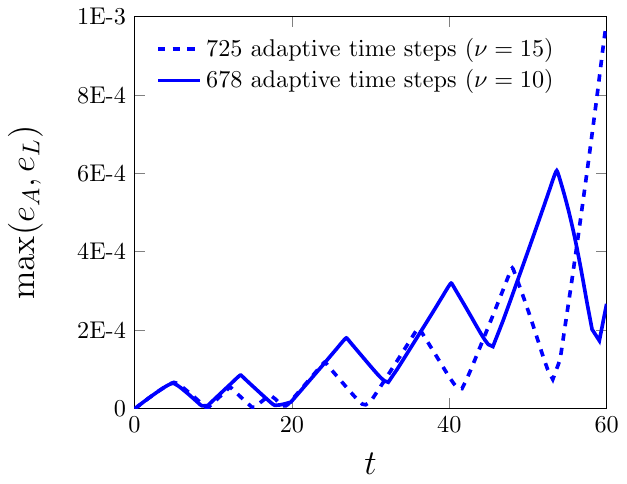}
\end{minipage}
\caption{\label{f:shear1VesErrors} {\em Left}: The error as a function
of time using a 150 uniform (dashed) and 381 adaptive (solid) time
steps for a single vesicle with $\nu=4$.  Even though the final errors
are comparable, the adaptive time stepping scheme does not incur a
large error near the start of the simulation.  {\em Right}: The error
as a function of time using 678 adaptive time steps for a single
vesicle with $\nu=10$ (solid).  The error happens to be decreasing near
the time horizon, and this explains why the final error is not tight to
the tolerance of 1E-3.  This is in contrast to the simulation with
$\nu=15$ (dashed) whose error is increasing near the time horizon,
therefore, allowing the final error to be much tighter to the
tolerance.} 
\end{figure}

We note two abnormalities in the results.  First, based on the final
errors, it appears that a fixed time step size is more efficient.  For
instance, when $\nu=4$, $m=150$ fixed time step sizes has a comparable
final error to $m=381$ adaptive time step sizes.  However, using a
fixed time step size, the error actually achieves a value of 4.04E-3
near the start of the simulation, and then decays to the final error
1.68E-3, whereas the adaptive time step error never exceeds 1E-3 (see
Fig.~\ref{f:shear1VesErrors}).  The errors for the other examples
behave similarly.  Second, for adaptive time stepping with $\nu=10$,
the final errors are much less than the requested tolerance.  This is a
consequence of the vesicle's error in area increasing and decreasing
throughout the simulation; near the time horizon, the error in area is
decreasing (see Fig.~\ref{f:shear1VesErrors}).  The result is that
there is an insignificant amount of time to commit enough error to
achieve the desired tolerance.

\subsection{Two vesicles in a shear flow}
\begin{figure}[htps]
\begin{minipage}{.4\textwidth}
\includegraphics{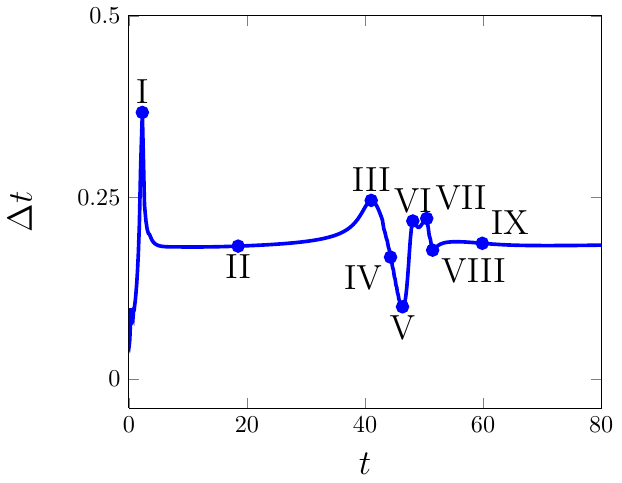}
\end{minipage}
\begin{minipage}{.58\textwidth}
\begin{tabular}{c|c|c}
\includegraphics[scale=0.32]{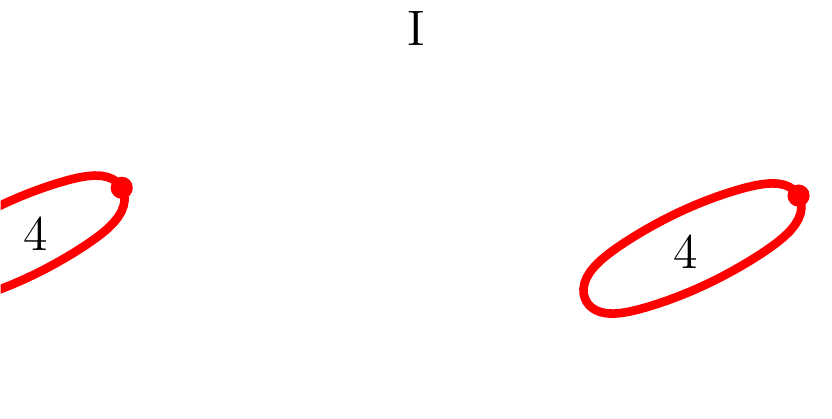} &
\includegraphics[scale=0.32]{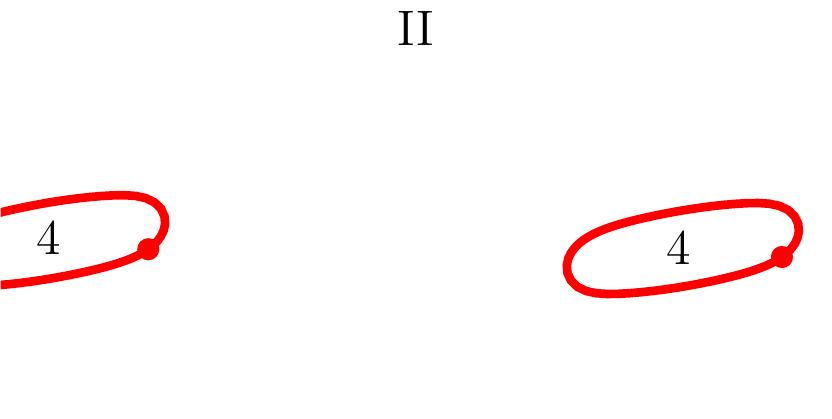} &
\includegraphics[scale=0.32]{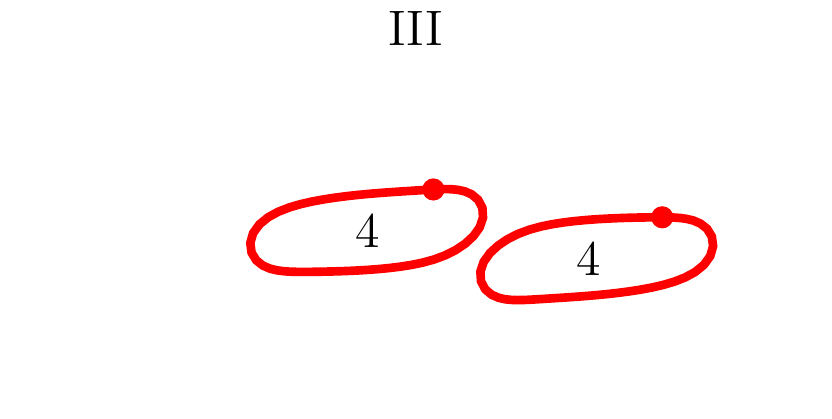} \\
\hline
\includegraphics[scale=0.32]{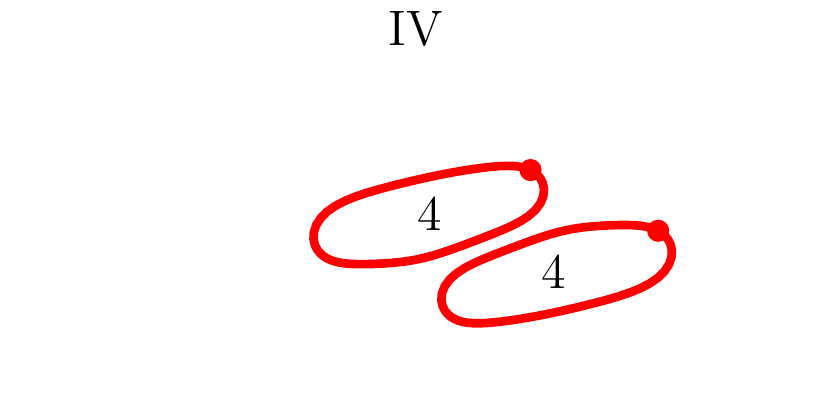} &
\includegraphics[scale=0.32]{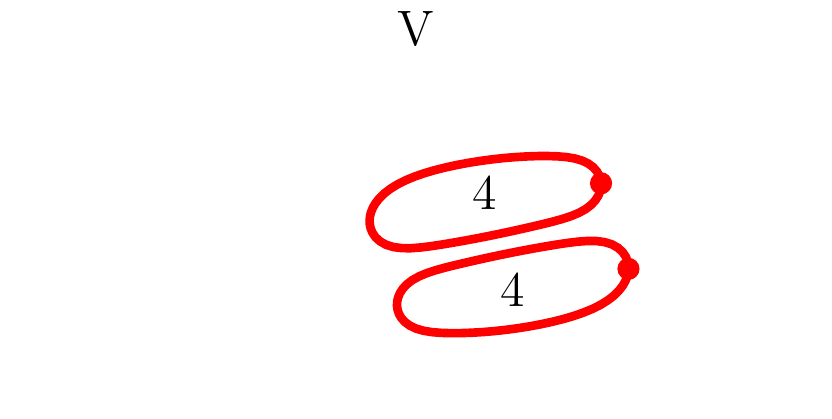} &
\includegraphics[scale=0.32]{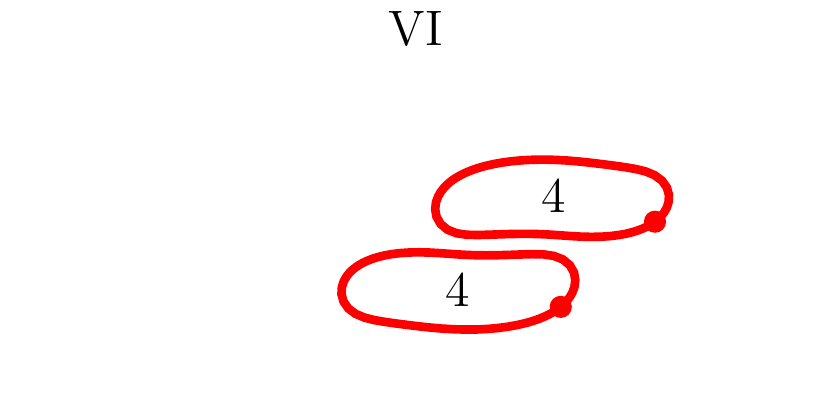} \\
\hline
\includegraphics[scale=0.32]{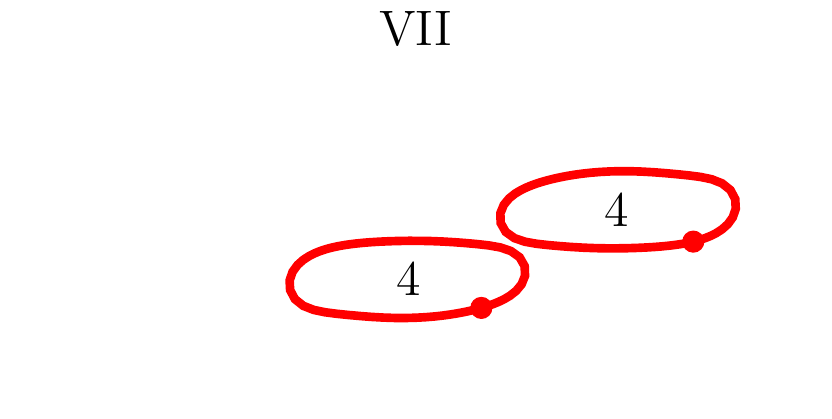} &
\includegraphics[scale=0.32]{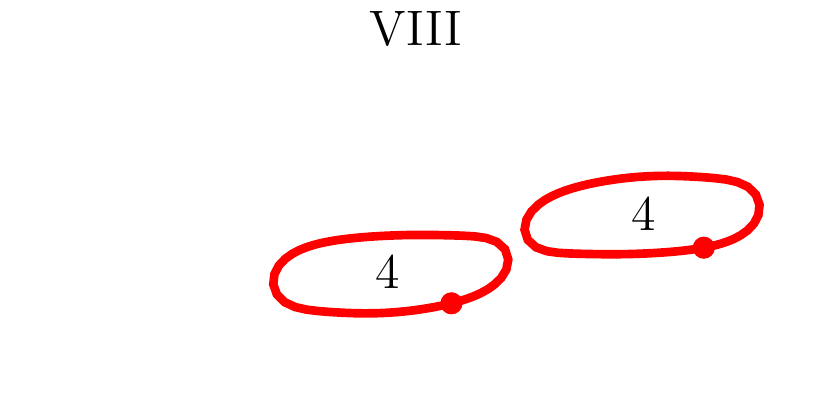} &
\includegraphics[scale=0.32]{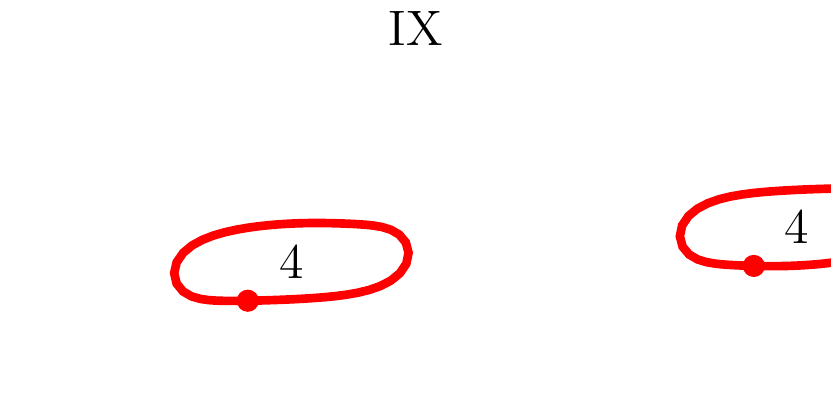} \\
\end{tabular}
\end{minipage}
\caption{\label{f:treadingANDtreading} Two tank-treading vesicles
($\nu=4$).  {\em Left}: The time step size taken to achieve the desired
tolerance 1.0E-2 at the time horizon.  The actual final error is
1.0E-2.  {\em Right}: The vesicle configuration at the indicated
labels.}
\begin{minipage}{.4\textwidth}
\includegraphics{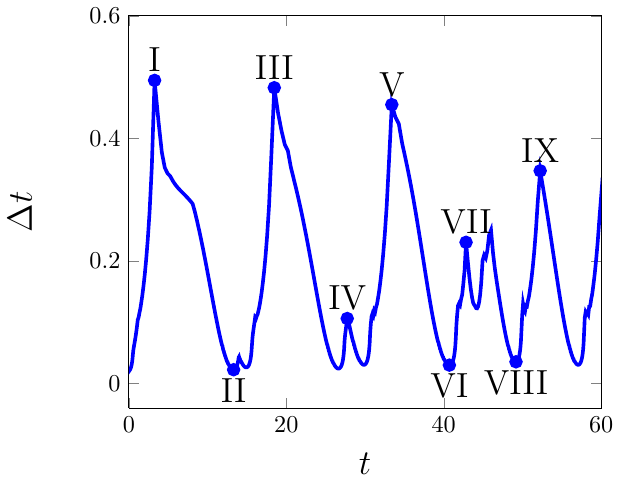}
\end{minipage}
\begin{minipage}{.58\textwidth}
\begin{tabular}{c|c|c}
\includegraphics[scale=0.32]{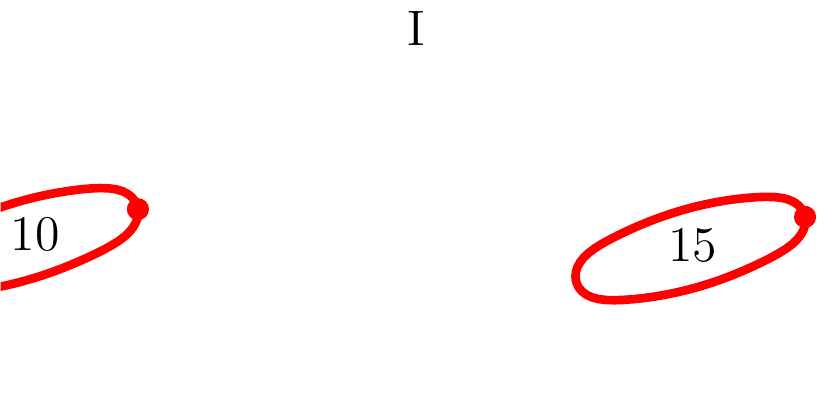} &
\includegraphics[scale=0.32]{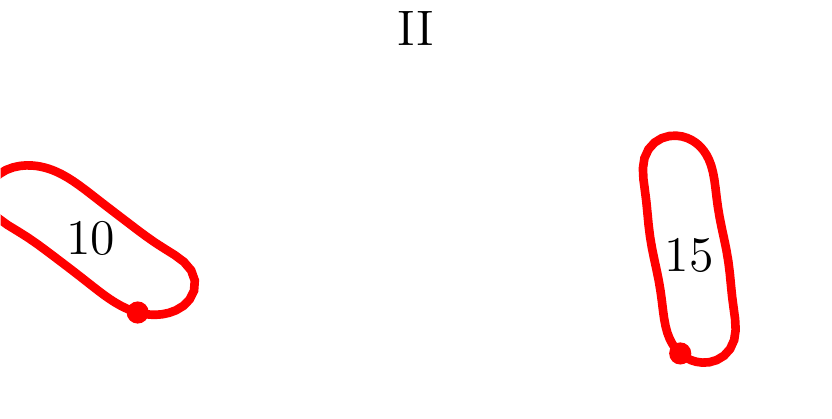} &
\includegraphics[scale=0.32]{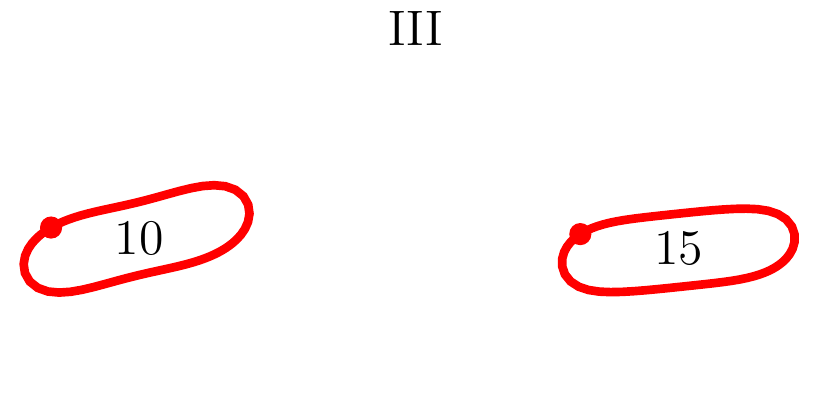} \\
\hline
\includegraphics[scale=0.32]{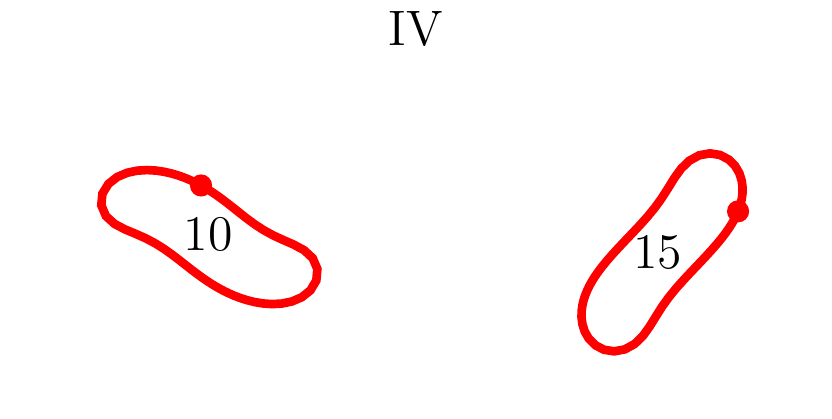} &
\includegraphics[scale=0.32]{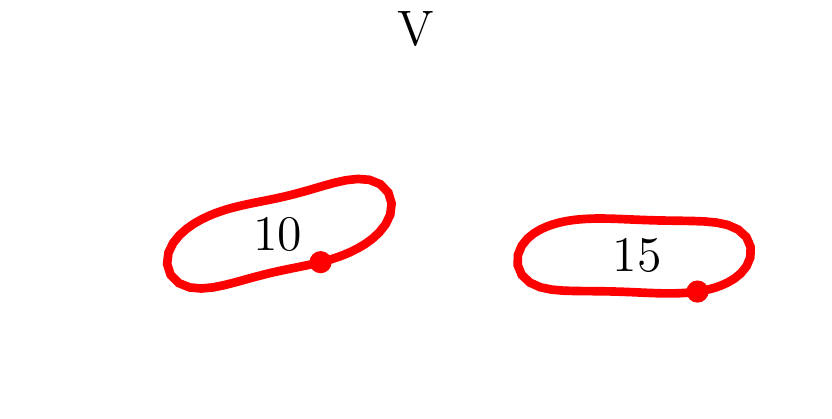} &
\includegraphics[scale=0.32]{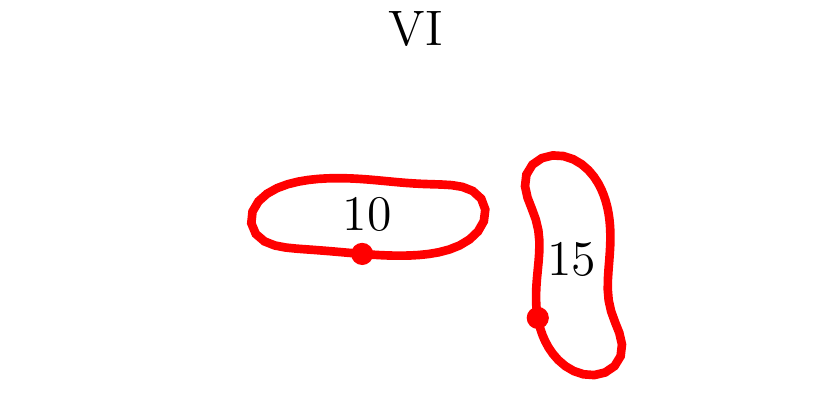} \\
\hline
\includegraphics[scale=0.32]{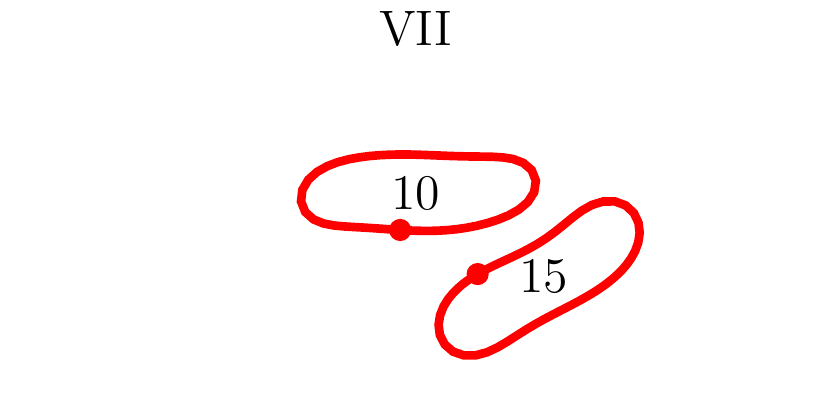} &
\includegraphics[scale=0.32]{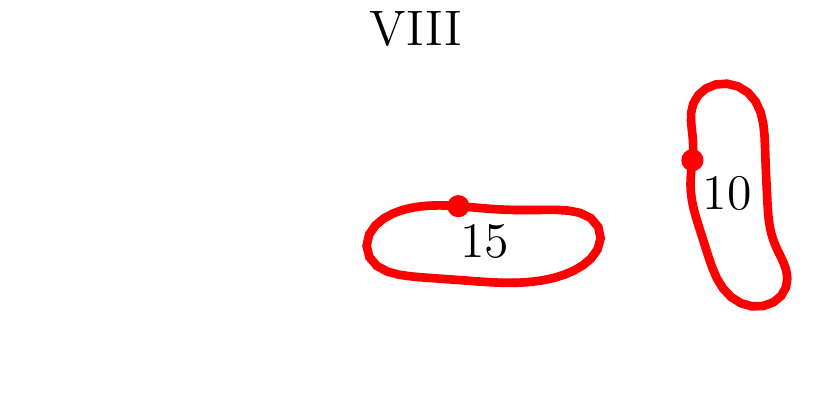} &
\includegraphics[scale=0.32]{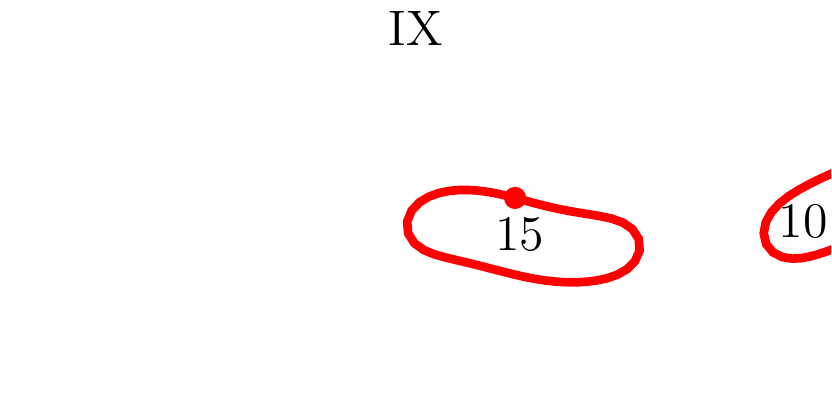} \\
\end{tabular}
\end{minipage}
\caption{\label{f:tumblingANDtumbling} Two tumbling vesicles ($\nu=10$
and $\nu=15$). {\em Left}: The time step size taken to achieve the
desired tolerance 1.0E-2 at the time horizon.  The actual final error
is 7.6E-3.  {\em Right}: The vesicle configuration at the indicated
labels.  The vesicles are labelled with their viscosity contrast.}
\begin{minipage}{.4\textwidth}
\includegraphics{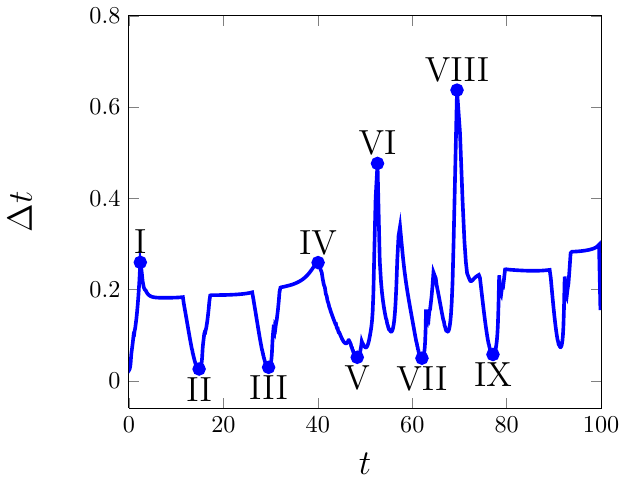}
\end{minipage}
\begin{minipage}{.58\textwidth}
\begin{tabular}{c|c|c}
\includegraphics[scale=0.32]{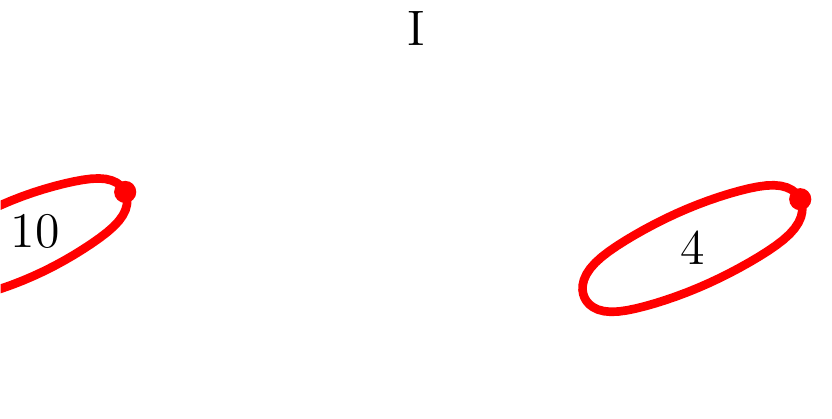} &
\includegraphics[scale=0.32]{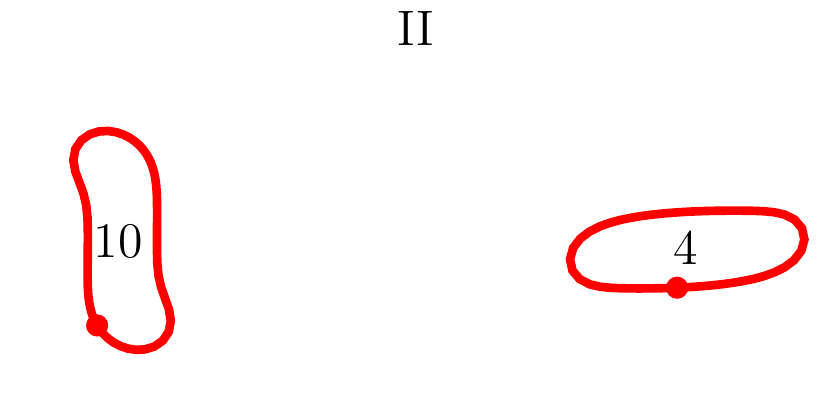} &
\includegraphics[scale=0.32]{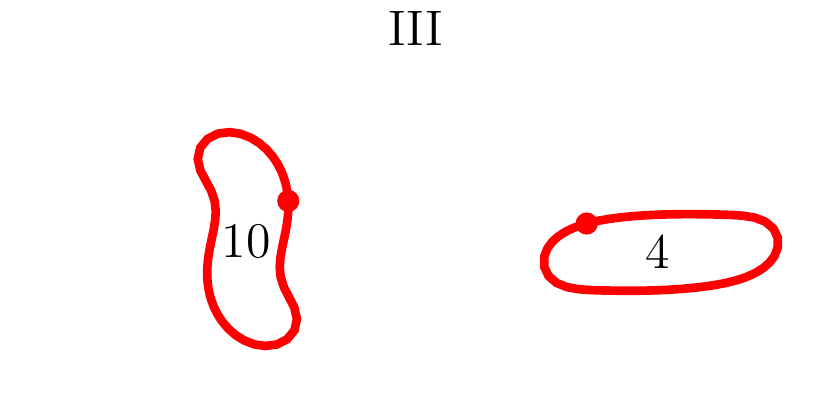} \\
\hline
\includegraphics[scale=0.32]{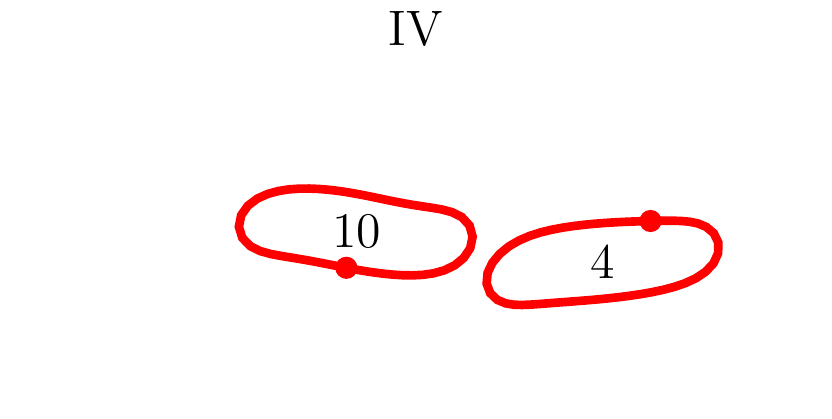} &
\includegraphics[scale=0.32]{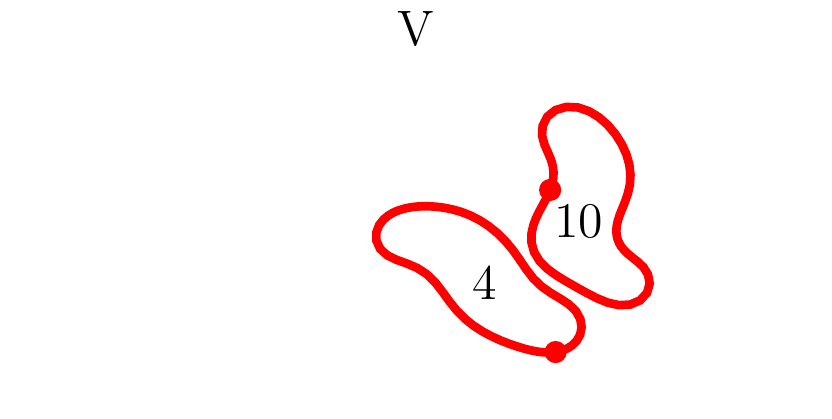} &
\includegraphics[scale=0.32]{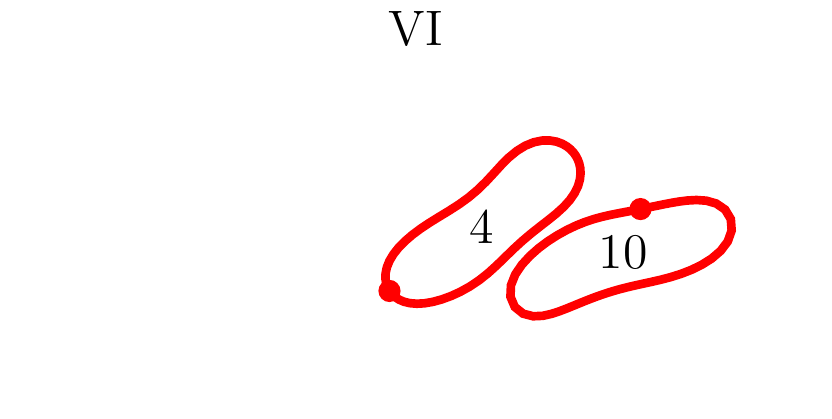} \\
\hline
\includegraphics[scale=0.32]{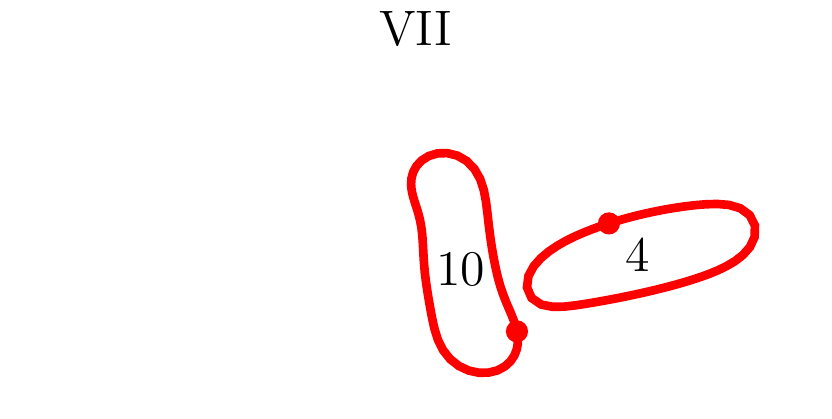} &
\includegraphics[scale=0.32]{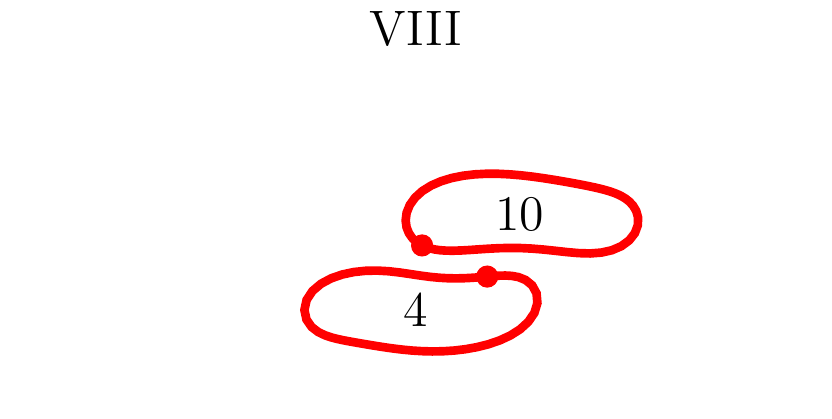} &
\includegraphics[scale=0.32]{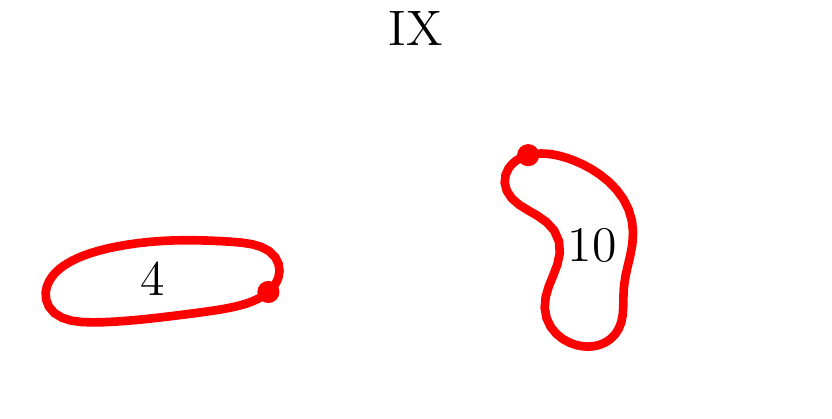} \\
\end{tabular}
\end{minipage}
\caption{\label{f:treadingANDtumbling} One tumbling vesicle ($\nu=10$)
and one tank-treading vesicle ($\nu=4$).  {\em Left}: The time step
size taken to achieve the desired tolerance 1.0E-2 at the time
horizon.  The actual final error is 1.0E-2.  {\em Right}: The vesicle
configuration at the indicated labels.  The vesicles are labelled with
their viscosity contrast.}
\end{figure}

We simulate two vesicles, each identical in shape to the last example,
where the left vesicle's center is slightly above the $x$-axis and the
right vesicle's center is on the $x$-axis.  We allow for larger spatial
errors by using $N=64$ points per vesicle, and the inter-vesicle
interactions are accelerated with the fast multipole
method~\cite{gre:rok1987}.  One SDC correction is used, and the
tolerance for the error in area and length at the time horizon is 1E-2.
The time step size and several vesicle configurations are illustrated in
Figs.~\ref{f:treadingANDtreading}--\ref{f:treadingANDtumbling} for
different viscosity contrasts.  Each vesicle is labelled with its
viscosity contrast and marked with a single tracker point to help
visualize the vesicle dynamics.

In Fig.~\ref{f:treadingANDtreading}, both the vesicles are tank-treading
($\nu=4$).  We observe that the time step size is nearly constant except
when the vesicles approach one another, and there, the time step size is
reduced.  The final error is 1.0E-2 which requires 545 accepted and 25
rejected time steps.  In Fig.~\ref{f:tumblingANDtumbling}, one vesicle
has a viscosity contrast $\nu=10$, while the other has a viscosity
contrast $\nu=15$, which results in the this vesicle tumbling faster.
The time step size is far less smooth, and the final error is 7.6E-3
which requires 975 accepted and 315 rejected time steps.  Finally, in
Fig.~\ref{f:treadingANDtumbling}, one vesicle is tumbling ($\nu = 10$),
while the other is tank-treading ($\nu = 4$).  We see that the time step
size is small when the tumbling vesicle has regions of high curvature,
and when the vesicles are close.  The final error is 1.0E-2 which
requires 752 accepted and 158 rejected time steps.

\section{Conclusion}
\label{s:Conclusion}
In this paper, we have extended our work on adaptive high-order time
integrators for vesicle suspensions~\cite{qua:bir2014c} to suspensions
that include viscosity contrast.  We have tested a second-order adaptive
time integrator on two different suspensions with varying viscosity
contrasts.  Future work includes spatial adaptivity and implementing the
presented results to three dimensions.

\section*{Acknowledgments}
This material is based upon work supported by AFOSR
grants FA9550-12-10484 and FA9550-11-10339; and NSF grants
CCF-1337393, OCI-1029022, and OCI-1047980; and by the U.S. Department
of Energy, Office of Science, Office of Advanced Scientific Computing
Research, Applied Mathematics program under Award Numbers
DE-SC0010518, DE-SC0009286, and DE-FG02-08ER2585.

\end{document}